\documentclass{amsart}
\usepackage{cases}
\usepackage{graphicx}
\usepackage{amssymb}
\usepackage{mathrsfs}

\newtheorem{theorem}{Theorem}[section]

\newtheorem{lemma}[theorem]{Lemma}

\numberwithin{equation}{section}

\newcommand{\bqn}{\begin{equation}}
\newcommand{\eqn}{\end{equation}}

\begin{document}
\title[]
{Sharp Bounds for Neuman Means in Terms of Geometric, Arithemtic and Quadratic Means}
\author{Zhi-Jun Guo}

\address{%
Zhi-Jun Guo, School of Mathematics and Computation Science, Hunan City University,
Yiyang 413000, China}
\email{492480045@qq.com}

\thanks{This research was supported by the Natural Science Foundation of China under Grants 11371125, 11171307 and 61374086, the Natural Science Foundation of Zhejiang Province under Grant LY13A010004, and the Natural Science Foundation of the Department of Education of Hunan Province under Grants 12C0577 and 13A013.}
\author{Yan Zhang}

\address{%
Yan Zhang, School of Mathematics and Computation Science, Hunan City University, Yiyang 413000, China}
\email{zhangyan080711zy@126.com}

\author{Yu-Ming Chu}

\address{%
Yu-Ming Chu (Corresponding author), School of Mathematics and Computation Science, Hunan City University, Yiyang 413000, China}
\email{chuyuming2005@126.com}

\author{Ying-Qing Song}

\address{%
Ying-Qing Song, School of Mathematics and Computation Science, Hunan City University, Yiyang 413000, China}
\email{1452225875@qq.com}


\begin{abstract}
In this paper, we find the greatest values  $\alpha_{1}$, $\alpha_{2}$, $\alpha_{3}$, $\alpha_{4}$, $\alpha_{5}$, $\alpha_{6}$, $\alpha_{7}$,
$\alpha_{8}$ and the least values $\beta_{1}$, $\beta_{2}$, $\beta_{3}$, $\beta_{4}$, $\beta_{5}$, $\beta_{6}$, $\beta_{7}$, $\beta_{8}$ such that the double inequalities
$$A^{\alpha_{1}}(a,b)G^{1-\alpha_{1}}(a,b)<N_{GA}(a,b)<A^{\beta_{1}}(a,b)G^{1-\beta_{1}}(a,b),$$
$$\frac{\alpha_{2}}{G(a,b)}+\frac{1-\alpha_{2}}{A(a,b)}<\frac{1}{N_{GA}(a,b)}<\frac{\beta_{2}}{G(a,b)}+\frac{1-\beta_{2}}{A(a,b)},$$
$$A^{\alpha_{3}}(a,b)G^{1-\alpha_{3}}(a,b)<N_{AG}(a,b)<A^{\beta_{3}}(a,b)G^{1-\beta_{3}}(a,b),$$
$$\frac{\alpha_{4}}{G(a,b)}+\frac{1-\alpha_{4}}{A(a,b)}<\frac{1}{N_{AG}(a,b)}<\frac{\beta_{4}}{G(a,b)}+\frac{1-\beta_{4}}{A(a,b)},$$
$$Q^{\alpha_{5}}(a,b)A^{1-\alpha_{5}}(a,b)<N_{AQ}(a,b)<Q^{\beta_{5}}(a,b)A^{1-\beta_{5}}(a,b),$$
$$\frac{\alpha_{6}}{A(a,b)}+\frac{1-\alpha_{6}}{Q(a,b)}<\frac{1}{N_{AQ}(a,b)}<\frac{\beta_{6}}{A(a,b)}+\frac{1-\beta_{6}}{Q(a,b)},$$
$$Q^{\alpha_{7}}(a,b)A^{1-\alpha_{7}}(a,b)<N_{QA}(a,b)<Q^{\beta_{7}}(a,b)A^{1-\beta_{7}}(a,b),$$
$$\frac{\alpha_{8}}{A(a,b)}+\frac{1-\alpha_{8}}{Q(a,b)}<\frac{1}{N_{QA}(a,b)}<\frac{\beta_{8}}{A(a,b)}+\frac{1-\beta_{8}}{Q(a,b)}$$
hold for all $a, b>0$ with $a\neq b$, where $G$, $A$ and $Q$ are respectively the geometric, arithmetic and quadratic means, and $N_{GA}$,
$N_{AG}$, $N_{AQ}$ and $N_{QA}$ are the Neuman means derived from the Schwab-Borchardt mean.

\footnotesize
\end{abstract}

\maketitle

\section{Introduction}
\bigskip

For $a,b>0$ with $a\neq b$, the Schwab-Borchardt mean $SB(a,b)$ [1-3] of $a$ and $b$ is
given by
\begin{equation*}
SB(a,b)=\begin{cases}
\frac{\sqrt{b^2-a^2}}{\cos^{-1}{(a/b)}}, &\quad  a<b, \\
\frac{\sqrt{a^2-b^2}}{\cosh^{-1}{(a/b)}}, &\quad a>b,
\end{cases}
\end{equation*}
where $\cos^{-1}(x)$ and $\cosh^{-1}(x)=\log(x+\sqrt{x^{2}-1})$ are the inverse cosine and inverse hyperbolic cosine functions, respectively. Recently, the Schwab-Borchardt mean has been the subject of intensive research. In particular, many remarkable inequalities for Schwab-Borchardt mean can be found in the literature [1-7]. Very recently, the Neuman mean $N(a,b)=(a+b^{2}/SB(a,b))/2$ derived from the Schwab-Borchardt was introduced and researched by Neuman in [8].

Let $N_{AG}(a,b)=N(A(a,b), G(a,b))$, $N_{GA}(a,b)=N(G(a,b), A(a,b))$, $N_{QA}(a,b)=N(Q(a,b), A(a,b))$ and $N_{AQ}(a,b)=N(A(a,b), Q(a,b))$ be the Neuman means, where $G(a,b)=\sqrt{ab}$, $A(a,b)=(a+b)/2$ and $Q(a,b)=\sqrt{(a^{2}+b^{2})/2}$ are the classical geometric, arithmetic and quadratic means of $a$ and $b$, respectively. Then Neuman [8] proved that the inequalities
$$G(a,b)<N_{AG}(a, b)<N_{GA}(a, b)<A(a,b)<N_{QA}(a, b)<N_{AQ}(a, b)<Q(a,b)$$
hold for all $a, b>0$ with $a\neq b$.

Let $a>b>0$ and $v=(a-b)/(a+b)\in (0, 1)$. Then we clearly see that
\begin{equation}
G(a,b)=A(a,b)\sqrt{1-v^{2}}, \quad Q(a,b)=A(a,b)\sqrt{1+v^{2}},
\end{equation}
and the following explicit formulas for $N_{AG}(a,b)$, $N_{GA}(a,b)$, $N_{QA}(a,b)$ and
$N_{AQ}(a, b)$ are given in [8]
\begin{equation}
N_{AG}(a,b)=\frac{1}{2}A(a,b)\left[1+(1-v^{2})\frac{\tanh^{-1}v}{v}\right],
\end{equation}
\begin{equation}
N_{GA}(a,b)=\frac{1}{2}A(a,b)\left[\sqrt{1-v^{2}}+\frac{\sin^{-1}v}{v}\right],
\end{equation}
\begin{equation}
N_{QA}(a,b)=\frac{1}{2}A(a,b)\left[\sqrt{1+v^{2}}+\frac{\sinh^{-1}v}{v}\right],
\end{equation}
\begin{equation}
N_{AQ}(a,b)=\frac{1}{2}A(a,b)\left[1+(1+v^{2})\frac{\tan^{-1}v}{v}\right],
\end{equation}
where $\tanh^{-1}(x)=\log[(1+x)/(1-x)]/2$, $\sin^{-1}(x)$, $\sinh^{-1}(x)=\log(x+\sqrt{1+x^{2}})$ and $\tan^{-1}(x)$ are the inverse hyperbolic tangent, inverse sine, inverse hyperbolic sine and inverse tangent functions, respectively.

In [8], Neuman also proved that the double inequalities
$$\alpha_{1} A(a,b)+(1-\alpha_{1})G(a,b)<N_{GA}(a,b)<\beta_{1} A(a,b)+(1-\beta_{1})G(a,b),$$
$$\alpha_{2} Q(a,b)+(1-\alpha_{2})A(a,b)<N_{AQ}(a,b)<\beta_{2} Q(a,b)+(1-\beta_{2})A(a,b),$$
$$\alpha_{3} A(a,b)+(1-\alpha_{3})G(a,b)<N_{AG}(a,b)<\beta_{3} A(a,b)+(1-\beta_{3})G(a,b),$$
$$\alpha_{4} Q(a,b)+(1-\alpha_{4})A(a,b)<N_{QA}(a,b)<\beta_{4} Q(a,b)+(1-\beta_{4})A(a,b)$$
hold for all $a, b>0$ with $a\neq b$ if and only if $\alpha_{1}\leq 2/3$, $\beta_{1}\geq \pi/4$,
$\alpha_{2}\leq 2/3$, $\beta_{2}\geq (\pi-2)/[4(\sqrt{2}-1)]=0.689\ldots$, $\alpha_{3}\leq 1/3$, $\beta_{3}\geq 1/2$,
$\alpha_{4}\leq 1/3$ and $\beta_{4}\geq [\log(1+\sqrt{2})+\sqrt{2}-2]/[2(\sqrt{2}-1)]=0.356\ldots$.

In [9], the authors presented the best possible parameters $\alpha_{1}, \alpha_{2}, \beta_{1}, \beta_{2}\in [0, 1/2]$ and
$\alpha_{3}, \alpha_{4}, \beta_{3}, \beta_{4}\in [1/2, 1]$ such that the double inequalities
$$G(\alpha_{1}a+(1-\alpha_{1})b, \alpha_{1}b+(1-\alpha_{1})a)<N_{AG}(a,b)<G(\beta_{1}a+(1-\beta_{1})b, \beta_{1}b+(1-\beta_{1})a),$$
$$G(\alpha_{2}a+(1-\alpha_{2})b, \alpha_{2}b+(1-\alpha_{2})a)<N_{GA}(a,b)<G(\beta_{2}a+(1-\beta_{2})b, \beta_{2}b+(1-\beta_{2})a),$$
$$Q(\alpha_{3}a+(1-\alpha_{3})b, \alpha_{3}b+(1-\alpha_{3})a)<N_{QA}(a,b)<Q(\beta_{3}a+(1-\beta_{3})b, \beta_{3}b+(1-\beta_{3})a),$$
$$Q(\alpha_{4}a+(1-\alpha_{4})b, \alpha_{4}b+(1-\alpha_{4})a)<N_{AQ}(a,b)<Q(\beta_{4}a+(1-\beta_{4})b, \beta_{4}b+(1-\beta_{4})a)$$
hold for all $a, b>0$ with $a\neq b$.

The main purpose of this paper is to find the greatest values  $\alpha_{1}$, $\alpha_{2}$, $\alpha_{3}$, $\alpha_{4}$, $\alpha_{5}$, $\alpha_{6}$, $\alpha_{7}$, $\alpha_{8}$ and the least values $\beta_{1}$, $\beta_{2}$, $\beta_{3}$, $\beta_{4}$, $\beta_{5}$, $\beta_{6}$, $\beta_{7}$, $\beta_{8}$ such that the double inequalities
$$A^{\alpha_{1}}(a,b)G^{1-\alpha_{1}}(a,b)<N_{GA}(a,b)<A^{\beta_{1}}(a,b)G^{1-\beta_{1}}(a,b),$$
$$\frac{\alpha_{2}}{G(a,b)}+\frac{1-\alpha_{2}}{A(a,b)}<\frac{1}{N_{GA}(a,b)}<\frac{\beta_{2}}{G(a,b)}+\frac{1-\beta_{2}}{A(a,b)},$$
$$A^{\alpha_{3}}(a,b)G^{1-\alpha_{3}}(a,b)<N_{AG}(a,b)<A^{\beta_{3}}(a,b)G^{1-\beta_{3}}(a,b),$$
$$\frac{\alpha_{4}}{G(a,b)}+\frac{1-\alpha_{4}}{A(a,b)}<\frac{1}{N_{AG}(a,b)}<\frac{\beta_{4}}{G(a,b)}+\frac{1-\beta_{4}}{A(a,b)},$$
$$Q^{\alpha_{5}}(a,b)A^{1-\alpha_{5}}(a,b)<N_{AQ}(a,b)<Q^{\beta_{5}}(a,b)A^{1-\beta_{5}}(a,b),$$
$$\frac{\alpha_{6}}{A(a,b)}+\frac{1-\alpha_{6}}{Q(a,b)}<\frac{1}{N_{AQ}(a,b)}<\frac{\beta_{6}}{A(a,b)}+\frac{1-\beta_{6}}{Q(a,b)},$$
$$Q^{\alpha_{7}}(a,b)A^{1-\alpha_{7}}(a,b)<N_{QA}(a,b)<Q^{\beta_{7}}(a,b)A^{1-\beta_{7}}(a,b),$$
$$\frac{\alpha_{8}}{A(a,b)}+\frac{1-\alpha_{8}}{Q(a,b)}<\frac{1}{N_{QA}(a,b)}<\frac{\beta_{8}}{A(a,b)}+\frac{1-\beta_{8}}{Q(a,b)}$$
hold for all $a, b>0$ with $a\neq b$ .

\bigskip
\section {Lemmas}
\bigskip

In order to prove our main results we need several lemmas, which we present in this section.

\setcounter{equation}{0}

\begin{lemma} (See [10, Theorem 1.25])  Let $-\infty<a<b<\infty$, $f, g: [a, b]\rightarrow (-\infty, \infty)$ be continuous on $[a,b]$ and differentiable on $(a, b)$, and $g^{\prime}(x)\neq 0$ on $(a,b)$. If $f'(x)/g'(x)$ is increasing (decreasing) on $(a,b)$, then so are
$$\frac{f(x)-f(a)}{g(x)-g(a)}, \quad\quad\ \frac{f(x)-f(b)}{g(x)-g(b)}.$$
If $f'(x)/g'(x)$ is strictly monotone, the the monotonicity in the conclusion is also strict.
\end{lemma}

\begin{lemma} (See [11, Lemma 1.1])  Suppose that the power series $f(x)=\sum_{n=0}^{\infty}a_{n}x^{n}$ and $g(x)=\sum_{n=0}^{\infty}b_{n}x^{n}$ have
the radius of convergence $r>0$ and $b_{n}>0$ for all $n\geq 0$. If the sequence $\{a_{n}/b_{n}\}$ is (strictly) increasing (decreasing) for all $n\geq 0$, then the function $f(x)/g(x)$ is also (strictly) increasing (decreasing) on $(0, r)$.
\end{lemma}

\begin{lemma} The function
\begin{equation}
f_{1}(x)=\frac{\log [\sin(2x)]-\log[2x+\sin(2x)]+\log 2}{\log(\cos x)}
\end{equation}
is strictly increasing from $(0, \pi/2)$ onto $(2/3, 1)$.
\end{lemma}
{\em Proof.} It follows from (3.1) that
\begin{equation}
f_{1}(0)=\frac{2}{3},
\end{equation}
\begin{equation}
f_{1}\left(\frac{\pi}{2}^{-}\right)=1.
\end{equation}

Let $g_{1}(x)=\log [\sin(2x)]-\log[2x+\sin(2x)]+\log 2$, $h_{1}(x)=\log(\cos x)$, $g_{2}(x)=\sin(2x)-2x\cos(2x)$ and $h_{2}(x)=[2x+\sin(2x)]\sin^{2}x$. Then simple computations lead to
\begin{equation}
g_{1}(0^{+})=h_{1}(0)=g_{2}(0)=h_{2}(0)=0,
\end{equation}
\begin{equation}
f_{1}(x)=\frac{g_{1}(x)}{h_{1}(x)}, \quad \frac{g'_{1}(x)}{h'_{1}(x)}=\frac{g_{2}(x)}{h_{2}(x)},
\end{equation}
\begin{equation}
\frac{g'_{2}(x)}{h'_{2}(x)}=\frac{1}{\frac{1}{2}+\frac{\sin(2x)}{2x}}.
\end{equation}

It is well known that the function $\sin x/x$ is strictly decreasing on $(0, \pi)$, hence equation (2.6) leads to the conclusion that the function
$g'_{2}(x)/h'_{2}(x)$ is strictly increasing on $(0, \pi/2)$. Therefore, Lemma 2.3 follows from Lemma 2.1 and (2.2)-(2.5) together with the monotonicity of $g'_{2}(x)/h'_{2}(x)$.

\bigskip
\begin{lemma} The function
\begin{equation}
f_{2}(x)=\frac{\log [2x+\sinh(2x)]-\log[\sinh(x)]-2\log 2}{\log[\cosh(x)]}
\end{equation}
is strictly increasing from $(0, \infty)$ onto $(1/3, 1)$.
\end{lemma}
{\em Proof.} It follows from (2.7) that
\begin{equation}
f_{2}(0^{+})=\frac{1}{3},
\end{equation}
\begin{equation}
\lim_{x\rightarrow \infty}f_{2}(x)=1.
\end{equation}

Let $g_{3}(x)=\log [2x+\sinh(2x)]-\log[\sinh(x)]-2\log 2$ and $h_{3}(x)=\log[\cosh(x)]$. Then simple computations lead to
\begin{equation}
f_{2}(x)=\frac{g_{3}(x)}{h_{3}(x)}, \quad g_{3}(0^{+})=h_{3}(0)=0,
\end{equation}
\begin{equation}
\frac{g'_{3}(x)}{h'_{3}(x)}=\frac{\sinh(4x)-4x\cosh(2x)+2\sinh(2x)-4x}{\sinh(4x)+4x\cosh(2x)-2\sinh(2x)-4x}
\end{equation}
\begin{equation*}
=\frac{\sum_{n=1}^{\infty}\frac{\left(2^{2n}-2n\right)2^{2n+2}}{(2n+1)!}x^{2n+1}}
{\sum_{n=1}^{\infty}\frac{\left(2^{2n}+2n\right)2^{2n+2}}{(2n+1)!}x^{2n+1}}
\end{equation*}
\begin{equation*}
=\frac{\sum_{n=0}^{\infty}\frac{\left(2^{2n+2}-2n-2\right)2^{2n+4}}{(2n+3)!}x^{2n}}
{\sum_{n=0}^{\infty}\frac{\left(2^{2n+2}+2n+2\right)2^{2n+4}}{(2n+3)!}x^{2n}}.
\end{equation*}
Let
\begin{equation}
a_{n}=\frac{\left(2^{2n+2}-2n-2\right)2^{2n+4}}{(2n+3)!}, \quad b_{n}=\frac{\left(2^{2n+2}+2n+2\right)2^{2n+4}}{(2n+3)!}.
\end{equation}
Then
\begin{equation}
b_{n}>0, \quad \frac{a_{n+1}}{b_{n+1}}-\frac{a_{n}}{b_{n}}=\frac{(3n+2)2^{2n+2}}{\left(2^{2n+3}+n+2\right)\left(2^{2n+1}+n+1\right)}>0
\end{equation}
for all $n\geq 0$.

It follows from Lemma 2.2 and (2.11)-(2.13) that the function $g'_{3}(x)/h'_{3}(x)$ is strictly increasing on $(0, \infty)$. Therefore, Lemma 2.4 follows from Lemma 2.1 and (2.8)-(2.10) together with the monotonicity of $g'_{3}(x)/h'_{3}(x)$.

\bigskip
\begin{lemma} The function
\begin{equation}
f_{3}(x)=\frac{2x-\sin(2x)}{(1-\cos x)[2x+\sin(2x)]}
\end{equation}
is strictly increasing from $(0, \pi/2)$ onto $(2/3, 1)$.
\end{lemma}
{\em Proof.} It follows from (2.14) that
\begin{equation}
f_{3}\left(0^{+}\right)=\frac{2}{3},
\end{equation}
\begin{equation}
f_{3}\left(\pi/2\right)=1.
\end{equation}

Let $g_{4}(x)=2x-\sin(2x)$ and $h_{4}(x)=(1-\cos x)[2x+\sin(2x)]$. Then simple computations lead to
\begin{equation}
f_{3}(x)=\frac{g_{4}(x)}{h_{4}(x)}, \quad g_{4}(0)=h_{4}(0)=0,
\end{equation}
\begin{equation*}
g'_{4}(x)=4\sin^{2}x,
\end{equation*}
\begin{equation*}
h'_{4}(x)=2\sin^{2}x\cos x-4\cos^{3}x+4\cos^{2}x+2x\sin x,
\end{equation*}
\begin{equation}
g'_{4}(0)=h'_{4}(0)=0,
\end{equation}
\begin{equation}
\frac{g''_{4}(x)}{h''_{4}(x)}=\frac{4}{9\cos x+\frac{x}{\sin x}-4},
\end{equation}
\begin{equation}
\left(9\cos x+\frac{x}{\sin x}\right)'=-8\sin x-\frac{[2x-\sin(2x)]\cos x}{2\sin^{2}x}<0
\end{equation}
for $x\in (0, \pi/2)$.

Therefore, Lemma 2.5 follows easily from Lemma 2.1 and (2.15)-(2.20).

\bigskip
\begin{lemma} The function
\begin{equation}
f_{4}(x)=\frac{\sinh(x)\cosh^{2}(x)-2\sinh(x)\cosh(x)+x\cosh(x)}{\sinh(x)\cosh^{2}(x)-\sinh(x)\cosh(x)+x\cosh(x)-x}
\end{equation}
is strictly increasing from $(0, \infty)$ onto $(1/3, 1)$.
\end{lemma}
{\em Proof.} It follows from (2.21) that
\begin{equation}
f_{4}(x)=\frac{\frac{1}{4}\sinh(3x)+\frac{1}{4}\sinh(x)-\sinh(2x)+x\cosh(x)}{\frac{1}{4}\sinh(3x)+\frac{1}{4}\sinh(x)-\frac{1}{2}\sinh(2x)+x\cosh(x)-x}
\end{equation}
\begin{equation*}
=\frac{\sum_{n=1}^{\infty}\frac{3^{2n+1}-2^{2n+3}+8n+5}{4[(2n+1)!]}x^{2n+1}}{\sum_{n=1}^{\infty}\frac{3^{2n+1}-2^{2n+2}+8n+5}{4[(2n+1)!]}x^{2n+1}}
\end{equation*}
\begin{equation*}
=\frac{\sum_{n=0}^{\infty}\frac{3^{2n+3}-2^{2n+5}+8n+13}{4[(2n+3)!]}x^{2n}}{\sum_{n=0}^{\infty}\frac{3^{2n+3}-2^{2n+4}+8n+13}{4[(2n+3)!]}x^{2n}}.
\end{equation*}
Let
\begin{equation}
a_{n}=\frac{3^{2n+3}-2^{2n+5}+8n+13}{4[(2n+3)!]}, \quad b_{n}=\frac{3^{2n+3}-2^{2n+4}+8n+13}{4[(2n+3)!]}.
\end{equation}
Then simple computations lead to
\begin{equation}
b_{n}>\frac{3^{2n+3}-2^{2n+4}}{4[(2n+3)!]}=\frac{2^{2n+3}\left[\left(\frac{3}{2}\right)^{2n+3}-2\right]}{4[(2n+3)!]}>0
\end{equation}
\begin{equation}
\frac{a_{n+1}}{b_{n+1}}-\frac{a_{n}}{b_{n}}=\frac{\left(135\times 3^{2n}-24n-31\right)2^{2n+4}}{\left(3^{2n+3}-2^{2n+4}+8n+13\right)\left(3^{2n+5}-2^{2n+6}+8n+21\right)}>0
\end{equation}
for all $n\geq 0$.

Note that
\begin{equation}
f_{4}(0^{+})=\frac{1}{3},  \quad \lim_{x\rightarrow \infty}f_{4}(x)=\lim_{n\rightarrow \infty}\frac{a_{n}}{b_{n}}=1.
\end{equation}

Therefore, Lemma 2.6 follows easily from Lemma 2.2 and (2.22)-(2.26).

\bigskip
\section{Main Results}
\bigskip
\begin{theorem} The double inequalities
\begin{equation}
A^{\alpha_{1}}(a,b)G^{1-\alpha_{1}}(a,b)<N_{GA}(a,b)<A^{\beta_{1}}(a,b)G^{1-\beta_{1}}(a,b),
\end{equation}
\begin{equation}
\frac{\alpha_{2}}{G(a,b)}+\frac{1-\alpha_{2}}{A(a,b)}<\frac{1}{N_{GA}(a,b)}<\frac{\beta_{2}}{G(a,b)}+\frac{1-\beta_{2}}{A(a,b)}
\end{equation}
holds for all $a, b>0$ with $a\neq b$ if and only if $\alpha_{1}\leq 2/3$, $\beta_{1}\geq 1$, $\alpha_{2}\leq 0$ and
$\beta_{2}\geq 1/3$.
\end{theorem}
{\em Proof.} We clearly see that inequalities (3.1) and (3.2) can be rewritten as
\begin{equation}
\left(\frac{A(a,b)}{G(a,b)}\right)^{\alpha_{1}}<\frac{N_{GA}(a,b)}{G(a,b)}<\left(\frac{A(a,b)}{G(a,b)}\right)^{\beta_{1}}
\end{equation}
and
\begin{equation}
1-\beta_{2}<\frac{\frac{1}{G(a,b)}-\frac{1}{N_{GA}(a,b)}}{\frac{1}{G(a,b)}-\frac{1}{A(a,b)}}<1-\alpha_{2},
\end{equation}
respectively.

Since both the geometric mean $G(a,b)$ and arithmetic mean $A(a,b)$ are symmetric and homogeneous of degree 1, without loss of generality, we assume that $a>b$. Let $v=(a-b)/(a+b)\in (0, 1)$. Then from (1.1) and (1.3) we know that inequalities (3.3) and (3.4) are equivalent to
\begin{equation}
\alpha_{1}<\frac{\log\left[\frac{1}{2}\left(1+\frac{\sin^{-1}(v)}{v\sqrt{1-v^{2}}}\right)\right]}{\log\frac{1}{\sqrt{1-v^{2}}}}<\beta_{1}
\end{equation}
and
\begin{equation}
1-\beta_{2}<\frac{\sin^{-1}v-v\sqrt{1-v^{2}}}{(1-\sqrt{1-v^{2}})(v\sqrt{1-v^{2}}+\sin^{-1}v)}<1-\alpha_{2},
\end{equation}
respectively.

Let $x=\sin^{-1}(v)$. Then $x\in (0, \pi/2)$,
\begin{equation}
\frac{\log\left[\frac{1}{2}\left(1+\frac{\sin^{-1}(v)}{v\sqrt{1-v^{2}}}\right)\right]}{\log\frac{1}{\sqrt{1-v^{2}}}}
=\frac{\log[\sin(2x)]-\log[2x+\sin(2x)]+\log 2}{\log(\cos x)},
\end{equation}
\begin{equation}
\frac{\sin^{-1}v-v\sqrt{1-v^{2}}}{(1-\sqrt{1-v^{2}})(v\sqrt{1-v^{2}}+\sin^{-1}v)}=\frac{2x-\sin(2x)}{(1-\cos x)[2x+\sin(2x)]}.
\end{equation}

Therefore, inequality (3.1) holds for all $a,b>0$ with $a\neq b$ follows from (3.5) and (3.7) together with Lemma 2.3,
and inequality (3.2) holds for all $a,b>0$ with $a\neq b$ follows from (3.6) and (3.8) together with Lemma 2.5.

\medskip
\begin{theorem} The double inequalities
\begin{equation}
A^{\alpha_{3}}(a,b)G^{1-\alpha_{3}}(a,b)<N_{AG}(a,b)<A^{\beta_{3}}(a,b)G^{1-\beta_{3}}(a,b),
\end{equation}
\begin{equation}
\frac{\alpha_{4}}{G(a,b)}+\frac{1-\alpha_{4}}{A(a,b)}<\frac{1}{N_{AG}(a,b)}<\frac{\beta_{4}}{G(a,b)}+\frac{1-\beta_{4}}{A(a,b)}
\end{equation}
hold for all $a, b>0$ with $a\neq b$ if and only if $\alpha_{3}\leq 1/3$, $\beta_{3}\geq 1$, $\alpha_{4}\leq 0$ and $\beta_{4}\geq 2/3$.
\end{theorem}

{\em Proof.} We clearly see that inequalities (3.9) and (3.10) can be rewritten as
\begin{equation}
\left(\frac{A(a,b)}{G(a,b)}\right)^{\alpha_{3}}<\frac{N_{AG}(a,b)}{G(a,b)}<\left(\frac{A(a,b)}{G(a,b)}\right)^{\beta_{3}}
\end{equation}
and
\begin{equation}
1-\beta_{4}<\frac{\frac{1}{G(a,b)}-\frac{1}{N_{AG}(a,b)}}{\frac{1}{G(a,b)}-\frac{1}{A(a,b)}}<1-\alpha_{4},
\end{equation}
respectively.

Without loss of generality, we assume that $a>b$. Let $v=(a-b)/(a+b)\in (0, 1)$. Then it follows from (1.1) and (1.2) that inequalities (3.11) and (3.12) are equivalent to
\begin{equation}
\alpha_{3}<\frac{\log\left[\frac{1}{\sqrt{1-v^{2}}}+\frac{\sqrt{1-v^{2}}}{v}\tanh^{-1}(v)\right]-\log 2}{\log\frac{1}{\sqrt{1-v^{2}}}}<\beta_{3}
\end{equation}
and
\begin{equation}
1-\beta_{4}<\frac{v+(1-v^{2})\tanh^{-1}(v)-2v\sqrt{1-v^{2}}}{(1-\sqrt{1-v^{2}})[v+(1-v^{2})\tanh^{-1}(v)]}<1-\beta_{4},
\end{equation}
respectively.
Let $x=\tanh^{-1}(v)\in (0, \infty)$. Then simple computations lead to
\begin{equation}
\frac{\log\left[\frac{1}{\sqrt{1-v^{2}}}+\frac{\sqrt{1-v^{2}}}{v}\tanh^{-1}(v)\right]-\log 2}{\log\frac{1}{\sqrt{1-v^{2}}}}
\end{equation}
\begin{equation*}
=\frac{\log[2x+\sinh(2x)]-\log[\sinh(x)]-2\log 2}{\log[\cosh(x)]}
\end{equation*}
and
\begin{equation}
\frac{v+(1-v^{2})\tanh^{-1}(v)-2v\sqrt{1-v^{2}}}{(1-\sqrt{1-v^{2}})[v+(1-v^{2})\tanh^{-1}(v)]}
\end{equation}
\begin{equation*}
=\frac{\sinh(x)\cosh^{2}(x)-2\sinh(x)\cosh(x)+x\cosh(x)}{\sinh(x)\cosh^{2}(x)-\sinh(x)\cosh(x)+x\cosh(x)-x}.
\end{equation*}

Therefore, inequality (3.9) holds for all $a, b>0$ with $a\neq b$ if and only if $\alpha_{3}\leq 1/3$ and $\beta_{3}\geq 1$ follows from (3.13) and (3.15) together with Lemma 2.4, and inequality (3.10) holds for all $a, b>0$ with $a\neq b$ if and only if $\alpha_{4}\leq 0$ and $\beta_{4}\geq 2/3$ follows from (3.14) and (3.16) together with Lemma 2.6.

\medskip
\begin{theorem} The double inequalities
\begin{equation}
Q^{\alpha_{5}}(a,b)A^{1-\alpha_{5}}(a,b)<N_{AQ}(a,b)<Q^{\beta_{5}}(a,b)A^{1-\beta_{5}}(a,b),
\end{equation}
\begin{equation}
\frac{\alpha_{6}}{A(a, b)}+\frac{1-\alpha_{6}}{Q(a,b)}<\frac{1}{N_{AQ}(a,b)}<\frac{\beta_{6}}{A(a, b)}+\frac{1-\beta_{6}}{Q(a,b)}
\end{equation}
hold for all $a, b>0$ with $a\neq b$ if and only if $\alpha_{5}\leq 2/3$, $\beta_{5}\geq 2\log(\pi+2)/\log 2-4=0.7244\ldots$, $\alpha_{6}\leq[6+2\sqrt{2}-(1+\sqrt{2})\pi]/(\pi+2)=0.2419\ldots$ and $\beta_{6}\geq 1/3$.
\end{theorem}

{\em Proof.} We clearly see that inequalities (3.17) and (3.18) can be rewritten as
\begin{equation}
\left(\frac{Q(a,b)}{A(a,b)}\right)^{\alpha_{5}}<\frac{N_{AQ}(a,b)}{A(a,b)}<\left(\frac{Q(a,b)}{A(a,b)}\right)^{\beta_{5}}
\end{equation}
and
\begin{equation}
1-\beta_{6}<\frac{\frac{1}{A(a,b)}-\frac{1}{N_{AQ}(a,b)}}{\frac{1}{A(a,b)}-\frac{1}{Q(a,b)}}<1-\alpha_{6},
\end{equation}
respectively.

Without loss of generality, we assume that $a>b$. Let $v=(a-b)/(a+b)\in (0, 1)$. Then from (1.1) and (1.5) we clearly see that inequalities (3.19) and (3.20) are equivalent to
\begin{equation}
\alpha_{5}<\frac{2\log(1+\frac{1+v^{2}}{v}\tan^{-1}(v))-2\log 2}{\log(1+v^{2})}<\beta_{5}
\end{equation}
and
\begin{equation}
1-\beta_{6}<\frac{\left[\left(1+v^{2}\right)\tan^{-1}(v)-v\right]\sqrt{1+v^{2}}}{\left[\left(1+v^{2}\right)\tan^{-1}(v)+v\right](\sqrt{1+v^{2}}-1)}<1-\alpha_{6},
\end{equation}
respectively.

Let $x=\tan^{-1}(v)$. Then $x\in (0, \pi/4)$,
\begin{equation}
\frac{2\log(1+\frac{1+v^{2}}{v}\tan^{-1}(v))-2\log 2}{\log(1+v^{2})}
\end{equation}
\begin{equation*}
=\frac{\log[\sin(2x)]-\log[2x+\sin(2x)]+\log 2}{\log(\cos x)}=f_{1}(x)
\end{equation*}
and
\begin{equation}
\frac{\left[\left(1+v^{2}\right)\tan^{-1}(v)-v\right]\sqrt{1+v^{2}}}{\left[\left(1+v^{2}\right)\tan^{-1}(v)+v\right](\sqrt{1+v^{2}}-1)}
\end{equation}
\begin{equation*}
\frac{2x-\sin(2x)}{(1-\cos x)[2x+\sin(2x)]}=f_{3}(x).
\end{equation*}
Note that
\begin{equation}
f_{1}\left(\frac{\pi}{4}\right)=\frac{2\log(\pi+2)}{\log 2}-4,
\end{equation}
\begin{equation}
f_{3}\left(\frac{\pi}{4}\right)=\frac{(2+\sqrt{2})(\pi-2)}{\pi+2}=1-\frac{6+2\sqrt{2}-(1+\sqrt{2})\pi}{\pi+2}.
\end{equation}

Therefore, inequality (3.17) holds for all $a, b>0$ with $a\neq b$ if and only if $\alpha_{5}\leq 2/3$ and $\beta_{5}\geq 2\log(\pi+2)/\log 2-4$ follows from (3.21), (3.23), (3.25) and Lemma 2.3, and inequality (3.18) holds for all $a, b>0$ with $a\neq b$ if and only if $\alpha_{6}\leq[6+2\sqrt{2}-(1+\sqrt{2})\pi]/(\pi+2)$ and $\beta_{6}\geq 1/3$ follows from (3.22), (3.24), (3.26) and Lemma 2.5.

\medskip
\begin{theorem} The double inequalities
\begin{equation}
Q^{\alpha_{7}}(a,b)A^{1-\alpha_{7}}(a,b)<N_{QA}(a,b)<Q^{\beta_{7}}(a,b)A^{1-\beta_{7}}(a,b),
\end{equation}
\begin{equation}
\frac{\alpha_{8}}{A(a, b)}+\frac{1-\alpha_{8}}{Q(a,b)}<\frac{1}{N_{QA}(a,b)}<\frac{\beta_{8}}{A(a, b)}+\frac{1-\beta_{8}}{Q(a,b)}
\end{equation}
hold for all $a, b>0$ with $a\neq b$ if and only if $\alpha_{7}\leq 1/3$, $\beta_{7}\geq 2\log[\sqrt{2}+\log(1+\sqrt{2})]/\log 2-2=0.3977\ldots$, $\alpha_{8}\leq[2+\sqrt{2}-(1+\sqrt{2})\log(1+\sqrt{2})]/[\sqrt{2}+\log(1+\sqrt{2})]=0.5603\ldots$ and $\beta_{8}\geq 2/3$.
\end{theorem}

{\em Proof.} We clearly see that inequalities (3.27) and (3.28) can be rewritten as
\begin{equation}
\left(\frac{Q(a,b)}{A(a,b)}\right)^{\alpha_{7}}<\frac{N_{QA}(a,b)}{A(a,b)}<\left(\frac{Q(a,b)}{A(a,b)}\right)^{\beta_{7}}
\end{equation}
and
\begin{equation}
1-\beta_{8}<\frac{\frac{1}{A(a,b)}-\frac{1}{N_{QA}(a,b)}}{\frac{1}{A(a,b)}-\frac{1}{Q(a,b)}}<1-\alpha_{8},
\end{equation}
respectively.

Without loss of generality, we assume that $a>b$. Let $v=(a-b)/(a+b)\in (0, 1)$. Then from (1.1) and (1.4) we clearly see that inequalities (3.29) and (3.30) are equivalent to
\begin{equation}
\alpha_{7}<\frac{2\log\left[\sqrt{1+v^{2}}+\frac{\sinh^{-1}(v)}{v}\right]-2\log 2}{\log\left(1+v^{2}\right)}<\beta_{7}
\end{equation}
and
\begin{equation}
1-\beta_{8}<\frac{\left[v\left(1+v^{2}\right)+\sqrt{1+v^{2}}\sinh^{-1}(v)\right]-2v\sqrt{1+v^{2}}}{(\sqrt{1+v^{2}}-1)\left[v\sqrt{1+v^{2}}+\sinh^{-1}(v)\right]}<1-\alpha_{8},
\end{equation}
respectively.

Let $x=\sinh^{-1}(v)$. Then $x\in (0, \log(1+\sqrt{2}))$,
\begin{equation}
\frac{2\log\left[\sqrt{1+v^{2}}+\frac{\sinh^{-1}(v)}{v}\right]-2\log 2}{\log\left(1+v^{2}\right)}
\end{equation}
\begin{equation*}
=\frac{\log[2x+\sinh(2x)]-\log[\sinh(x)]-2\log 2}{\log[\cosh(x)]}=f_{2}(x),
\end{equation*}
\begin{equation}
\frac{\left[v\left(1+v^{2}\right)+\sqrt{1+v^{2}}\sinh^{-1}(v)\right]-2v\sqrt{1+v^{2}}}{(\sqrt{1+v^{2}}-1)\left[v\sqrt{1+v^{2}}+\sinh^{-1}(v)\right]}
\end{equation}
\begin{equation*}
=\frac{\sinh(x)\cosh^{2}(x)-2\sinh(x)\cosh(x)+x\cosh(x)}{\sinh(x)\cosh^{2}(x)-\sinh(x)\cosh(x)+x\cosh(x)-x}=f_{4}(x).
\end{equation*}
Note that
\begin{equation}
f_{2}[\log(1+\sqrt{2})]=\frac{2\log[\sqrt{2}+\log(1+\sqrt{2})]}{\log 2}-2,
\end{equation}
\begin{equation}
f_{4}[\log(1+\sqrt{2})]=\frac{(2+\sqrt{2})\log(1+\sqrt{2})-2}{\sqrt{2}+\log(1+\sqrt{2})}
\end{equation}
\begin{equation*}
=1-\frac{2+\sqrt{2}-(1+\sqrt{2})\log(1+\sqrt{2})}{\sqrt{2}+\log(1+\sqrt{2})}.
\end{equation*}

Therefore, inequality (3.27) holds for all $a, b>0$ with $a\neq b$ if and only if $\alpha_{7}\leq 1/3$ and $\beta_{7}\geq 2\log[\sqrt{2}+\log(1+\sqrt{2})]/\log 2-2$ follows from (3.31), (3.33), (3.35) and Lemma 2.4, and inequality (3.28) holds
for all $a, b>0$ with $a\neq b$ if and only if $\alpha_{8}\leq[2+\sqrt{2}-(1+\sqrt{2})\log(1+\sqrt{2})]/[\sqrt{2}+\log(1+\sqrt{2})]$ and $\beta_{8}\geq 2/3$ follows from (3.32), (3.34), (3.36) and Lemma 2.6.

\medskip

\noindent{\bf Conflict of Interests}

\noindent{The authors declare that there is no conflict of interests regarding the publication of this paper.}

\end{document}